\input amstex\documentstyle{amsppt}  
\pagewidth{12.5cm}\pageheight{19cm}\magnification\magstep1
\topmatter
\title On $C$-small conjugacy classes in a reductive group\endtitle
\author G. Lusztig\endauthor
\address{Department of Mathematics, M.I.T., Cambridge, MA 02139}\endaddress
\thanks{Supported in part by the National Science Foundation}\endthanks
\endtopmatter   
\document
\define\Irr{\text{\rm Irr}}
\define\sneq{\subsetneqq}

\define\dw{\dot w}

\define\dv{\dot v}

\define\dz{\dot z}

\define\ul{\un l}

\define\uWW{\un\WW}

\define\mpb{\medpagebreak}

\define\da{\dagger}

\define\sqc{\sqcup}

\define\qua{\quad}

\define\lb{\linebreak}

\define\op{\oplus}

\define\part{\partial}
\define\em{\emptyset}

\define\m{\mapsto}
\define\do{\dots}

\define\sm{\smallmatrix}
\define\esm{\endsmallmatrix}
\define\sub{\subset}    
\define\bxt{\boxtimes}
\define\T{\times}
\define\ti{\tilde}
\define\nl{\newline}
\redefine\i{^{-1}}

\define\un{\underline}

\define\ot{\otimes}
\define\bbq{\bar{\QQ}_l}

\define\ind{\text{\rm ind}}

\define\sg{\text{\rm sgn}}
\define\tr{\text{\rm tr}}

\define\a{\alpha}

\redefine\c{\chi}
\define\g{\gamma}
\redefine\d{\delta}

\define\p{\pi}
\define\ph{\phi}
\define\ps{\psi}
\define\r{\rho}
\define\s{\sigma}
\redefine\t{\tau}

\define\k{\kappa}
\redefine\l{\lambda}
\define\z{\zeta}

\redefine\D{\Delta}
\define\Om{\Omega}

\define\Ph{\Phi}

\define\kk{\bold k}

\define\nn{\bold n}

\define\FF{\bold F}

\define\NN{\bold N}

\define\QQ{\bold Q}

\define\WW{\bold W}
\define\ZZ{\bold Z}

\define\cb{\Cal B}

\define\cf{\Cal F}

\define\co{\Cal O}

\define\cs{\Cal S}
\define\ct{\Cal T}

\define\cv{\Cal V}
\define\cw{\Cal W}
\define\cz{\Cal Z}

\define\fB{\frak B}

\define\fK{\frak K}

\define\tp{\ti p}

\define\sha{\sharp}

\define\bp{\bar p}

\define\bul{\bullet}

\define\DL{DL}
\define\GP{GP}
\define\GH{C}
\define\OR{L1}
\define\IC{L2}
\define\CSI{L3}
\define\CSIII{L4}
\define\WEUN{L5}
\define\SPE{L6}
\define\SPS{Sp}

\head Introduction\endhead
\subhead 0.1\endsubhead
Let $G$ be a connected reductive algebraic group over an algebraically closed field $\kk$ of characteristic $p$.
Let $\WW$ be the Weyl group of $G$. Let $\Om_w$ be the double coset in $G$ (with respect to a Borel subgroup 
$B^*$) corresponding to an element $w\in\WW$ which has minimal length in its conjugacy class $C$ in $\WW$ and has
no eigenvalue $1$ in the reflection representation of $\WW$. Let $Z_G$ be the centre of $G$.
Let $\g$ be a conjugacy class in $G$ such that $\Om\cap\g\ne\em$.
Since $\Om_w\cap\g$ is $B^*/Z_G$-stable for the conjugation action of $B^*/Z_G$ and this action
has finite isotropy groups by \cite{\WEUN, 5.2}, we see that $\dim(\Om_w\cap\g)\ge\dim(B^*/Z_G)$. 
As in \cite{\WEUN} we say that $\g$ is $C$-small if the previous inequality is an equality. (This condition 
depends only on $C$, not on $w$.) 

In the remainder of this subsection we assume that 

(i) $p$ is $0$ or a good prime for $G$
\nl
and that $G$ is almost
simple. In \cite{\WEUN} we have shown that for any $C$ as above there is a unique unipotent class $\g_C$ in $G$ 
which is $C$-small. In this paper we investigate the existence of $C$-small semisimple classes in $G$. We show 
that such a class $\g'$ exists in almost all cases. (There is exactly one exception to this property: it arises 
in type $E_8$ for a unique $C$.) Let $\r_{\g_C}$ be the Springer representation of $\WW$ associated to $\g_C$ and
to the local system $\bbq$ on $\g_C$. We show that $\r_{\g_C}$ is surprisingly connected to $\g'$ above (again 
with the unique exception above) as follows: $\r_{\g_C}$ is obtained by "$j$-induction" (see 0.3) from the sign 
representation of a reflection subgroup of $\WW$, namely the Weyl group of the connected centralizer of an 
element of $\g'$. We will also show that the representation $\r_{\g_C}$ depends only on the Weyl group $\WW$, not
on the underlying root system.

\subhead 0.2\endsubhead
Here is some notation that we use in this paper. Let $\cb$ the variety of Borel subgroups of $G$.
Let $\ul:\WW@>>>\NN$ be the standard length function. Let $S=\{s\in\WW;\ul(s)=1\}$. For 
each $w\in\WW$ let $\co_w$ be the corresponding $G$-orbit in $\cb\T\cb$. Let $\uWW_{el}$ be the set of elliptic 
conjugacy classes in $\WW$ (see \cite{\WEUN, 0.2}.) For $C\in\uWW_{el}$ let 
$d_C=\min_{w\in C}\ul(w)$ and let $C_{min}=\{w\in C;l(w)=d_C\}$. For any conjugacy class $\g$ in $G$ and any
$w\in\WW$ we set $\fB_w=\{(g,B)\in G\T\cb;(B,gBg\i)\in\co_w\}$, $\fB^\g_w=\{(g,B)\in\fB_w;g\in\g\}$. 

The cardinal of a finite set $X$ is denoted by $|X|$ or by $\sha(X)$. For any $g\in G$ let $Z(g)$ be the
centralizer of $g$ in $G$. Let $\nu_G$ be the number of positive roots of $G$. For an integer $\s$ we define 
$\k_\s\in\{0,1\}$ by $\s=\k_\s\mod2$.

Let $C\in\uWW_{el}$. In \cite{\WEUN, 5.5, 5.7(iii)} it is shown that if $\fB^\g_w\ne\em$ for some/any 
$w\in C_{min}$ then $\dim\fB^\g_w\ge\dim(G/Z_G)$ and $\dim\g\ge\dim(G/Z_G)-d_C$. (Here the equivalence of 
"some/any" follows from \cite{\WEUN, 5.2(a)}.) Following \cite{\WEUN, 5.5} we say that $\g$ is {\it $C$-small} if
for some/any $w\in C_{min}$ we have $\fB^\g_w\ne\em$ and the equivalent conditions $\dim\fB^\g_w=\dim(G/Z_G)$,
$\dim\g=\dim(G/Z_G)-d_C$ are satisfied (for the equivalence see \cite{\WEUN, 7.7(iv)}).

\subhead 0.3\endsubhead
Let $W$ be a Weyl group. Let $\sg$ be the sign representation of $W$ and let $R_W$ be the reflection 
representation of $W$. Let $\Irr(W)$ be the set of 
(isomorphism classes of) irreducible representations of $W$. For $E\in\Irr(W)$ let $b_E$ be the smallest integer 
$\ge0$ such that the multiplicity of $E$ in the $b_E$-th symmetric power of $R_W$ is $\ge1$. We write 
$E\in\Irr(W)^\da$ if this multiplicity is $1$. Let $W'$ be a subgroup of $W$ generated by reflections. Let 
$E'\in\Irr(W')^\da$. There is a unique $E\in\Irr(W)$ such that $E$ appears in $\ind_{W'}^W(E')$ and $b_E=b_{E'}$.
(See \cite{\GP, 5.2.6}.) We have $E\in\Irr(W)^\da$. We set $E=j_{W'}^W(E')$. The process $j_{W'}^W()$ is called
$j$-induction.

\subhead 0.4\endsubhead
For any unipotent class $\g$ in $G$ let $\r_\g\in\Irr(\WW)$ be the Springer representation of $\WW$ associated to
$\g$ and the local system $\bbq$ on $\g$. (We use the conventions of \cite{\IC}.) For any $C\in\uWW_{el}$ let 
$\g_C$ be the unique $C$-small unipotent class of $G$, see \cite{\WEUN}; thus $\r_{\g_C}$ is well defined.

\subhead 0.5\endsubhead
Let $B^*$ be a Borel subgroup of $G$ and let $\ct$ be a maximal torus of $B^*$. Let $N_G(\ct)$ be the normalizer 
of $\ct$ in $G$ and let $\cw=N_G(\ct)/\ct$. For any $z\in\cw$ let $\dz$ be a representative of $z$ in $N_G(\ct)$.
We identify $\WW=\cw$ as follows: to $z\in\cw$ corresponds the element $w\in\WW$ such that 
$(B^*,\dz B^*\dz\i)\in\co_w$. For any $s\in S$ let $\a_s:\ct@>>>\kk^*$ be the simple root defined by $s$. In the 
remainder of this subsection we assume that $G$ is almost simple, simply connected and that 0.1(i) holds.
Let $\a_0:\ct@>>>\kk^*$ be 
the unique root such that for any $s\in S$, $\a_0\a_s\i:\ct@>>>\kk^*$ is not a root. Let 
$\D=\{\a_s;s\in S\}\sqc\{\a_0\}$. For any $K\sneq\D$ let $\cw_K$ be the subgroup of $\cw$ generated by the 
reflections with respect to roots in $\D$. Let $G_K$ be the subgroup of $G$ generated by $\ct$ and by the root 
subgroups attached to roots such that the corresponding reflection in $\cw$ is in $\cw_K$. Note that $G_K$ is a 
connected reductive subgroup of $G$ (a "Borel-de Siebenthal subgroup"). Now $B^*\cap G_K$ is a Borel subgroup of 
$G_K$ and $\ct$ is a maximal torus of $B^*\cap G_K$. Note that $\cw_K=N_{G_K}(\ct)/\ct=\{w\in\cw;\dw\in G_K\}$ 
may be identified (using $B^*\cap G_K$, $\ct$) with the Weyl group $\WW_K$ of $G_K$ in the same way as $\cw$ is 
identified with $\WW$ (using $B^*,\ct$). In particular $\WW_K$ appears as a subgroup of $\WW$. Let $\cs_K$ be the
set of semisimple conjugacy classes $\g$ in $G$ such that for some $\z\in\g\cap\ct$ we have $G_K=Z(\z)$. Note 
that $\cs_K\ne\em$ and any semisimple class in $G$ belongs to $\cs_K$ for some $K\sneq\D$. The following is our
main result.

\proclaim{Theorem 0.6} Assume that $G$ is almost simple, simply connected and that 0.1(i) holds. Let 
$C\in\uWW_{el}$. With the single exception when $G$ is of type $E_8$ and for any $w\in C$, the characteristic 
polynomial of $w:R_\WW@>>>R_\WW$ is $(X+1)(X^2+1)^2(X^3+1)$, there exists $K\sneq\D$ such that 

(i) for any $\g\in\cs_K$, $\g$ is a $C$-small semisimple class;

(ii) $\r_{\g_C}=j_{\WW_K}^\WW(\sg)$.
\endproclaim
In the case where $G$ is of type $E_8$ and $C$ is the class specified in the theorem, there is no $K\sneq\D$ for 
which (i) holds and there is no $K\sneq\D$ for which (ii) holds. On the other hand in this case we have 
$\r_{\g_C}=j_{\WW_K}^\WW(\sg\ot r)$ where $G_K$ is of type $D_5+A_3$ and $r$ is the irreducible representation of
$\WW_K$ on which the $D_5$-factor acts as the reflection representation and the $A_3$-factor acts trivially. 
Also, if $\g\in\cs_K$, $\z\in\g\cap\ct$, $Z(\z)=K$ and $u$ is a unipotent element of $G_K$ which is in a minimal 
unipotent class $\ne1$ of the $D_5$-factor, then the $G$-conjugacy class $\g'$ of $\z u$ is $C$-small; although 
$\g'$ is not semisimple, it is as close as possible to being semisimple.
 
In the case where $G$ is of type $A$ the theorem is immediate: $C$ must be the conjugacy class of a Coxeter 
element and we can take $K=\em$. The proof of the theorem in the case where $G$ is of classical type other than 
$A$ is given in \S1, \S3. When $G$ is of exceptional type the proof of the theorem is given in 2.4, 3.5 (using a 
reduction to a computer calculation, see 2.2.)

\subhead 0.7\endsubhead
Assume that $G$ is almost simple, simply connected and that 0.1(i) holds.
Let $\g$ be any $C$-small conjugacy class in $G$. Let $\z$ (resp. $u$) be a semisimple (resp. unipotent) element
of $G$ such that $\z u=u\z\in\g$. Let $K\sneq\D$ be such that the conjugacy class of $\z$ belongs to $\cs_K$. We
assume as we may that $\z\in\g\cap\ct$ and $G_K=Z(\z)$. Let $\r_u$ be the Springer representation of $\WW_K$
associated to the conjugacy class of $u$ in $G_K$. We conjecture that 
$$\r_{\g_C}=j_{\WW_K}^\WW(\r_u).$$
This is supported by Theorem 0.6.

\subhead 0.8\endsubhead
I thank Gongqin Li for her help with programming in GAP3.

\head 1. Classical groups\endhead
\subhead 1.1\endsubhead
Let $t_1,t_2,\do,t_m$ be commuting indeterminates. Let $M$ be the $m\T m$ matrix whose $j$-th row is 
$t_1^j-t_1^{-j},t_2^j-t_2^{-j},\do,t_m^j-t_m^{-j}$, $j\in[1,m]$. Let $M'$ be the $(m+2)\T(m+2)$ matrix whose 
$j$-th row is $1,(-1)^j,t_1^j+t_1^{-j},t_2^j+t_2^{-j},\do,t_m^j+t_m^{-j}$, $j\in[0,m+1]$. The proof of (a),(b) 
below is left to the reader.

(a) {\it $\det(M)$ is equal to $\pm\prod_i\prod_{i<j}(t_i-t_j)\prod_{i\le j}(t_i-t_j\i)$ times a monomial in the 
$t_i$;}

(b) {\it $\det(M')$ is equal to $\pm2\prod_i(t_i-t_i\i)\prod_{i<j}(t_i-t_j)\prod_{i\le j}(t_i-t_j\i)$ times a 
monomial in the $t_i$.}

\subhead 1.2\endsubhead
Let $V$ be a $\kk$-vector space of finite dimension $\nn\ge3$. We set $\k=\k_\nn$ so that $\nn=2n+\k$ with
$n\in\NN$. Assume that $V$ has a fixed bilinear form $(,):V\T V@>>>\kk$ and a fixed quadratic form $Q:V@>>>\kk$ 
such that (i) or (ii) below holds:

(i) $Q=0$, $(x,x)=0$ for all $x\in V$, $(,)$ is nondegenerate;

(ii) $Q\ne0$, $(x,y)=Q(x+y)-Q(x)-Q(y)$ for $x,y\in V$, $p\ne2$, $(,)$ is nondegenerate.
\nl
An element $g\in GL(V)$ is said to be an isometry if $(gx,gy)=(x,y)=0$ for all $x,y\in V$ (hence $Q(gx)=Q(x)$ for
all $x\in V$). Let $Is(V)$ be the group of all isometries of $V$ (a closed subgroup of $GL(V)$). In this section 
we assume that $G$ is the identity component of $Is(V)$. Let $\cf$ be the set of all sequences 
$V_*=(0=V_0\sub V_1\sub V_2\sub\do\sub V_\nn=V)$ of subspaces of $V$ such that $\dim V_i=i$ for $i\in[0,\nn]$, 
$Q|_{V_i}=0$ and $\{x\in V;(x,V_i)=0\}=V_{\nn-i}$ for all $i\in[0,n]$. Now $Is(V)$ acts naturally (transitively) 
on $\cf$. 

\subhead 1.3\endsubhead
Assume that $Q=0$ so that $\nn=2n$.
Let $p_*=(p_1\ge p_2\ge\do\ge p_\s)$ be a sequence of integers $\ge1$ such that $p_1+p_2+\do+p_\s=n$.
For any $i\ge1$ we set $\bp_i=\sha(t\in[1,\s];p_t\ge i)$ so that $\bp_1\ge\bp_2\ge\do$ and $\sum_i\bp_i=n$. Let
$k=p_1$. We have $\bp_k\ge1$, $\bp_{k+1}=0$. We can find subspaces $\cv_i,\cv'_i$ ($i\in[1,k]$) of $V$ such that 

$V=\cv_1\op\cv'_1\op\cv_2\op\cv'_2\op\do\op\cv_k\op\cv'_k$;

$\dim\cv_i=\dim\cv'_i=\bp_i$ for $i\in[1,k]$;

$(,)$ is zero on $\cv_i,\cv'_i$ for $i\in[1,k]$;

$(\cv_i\op\cv'_i,\cv_j\op\cv'_j)=0$ for all $i\ne j$.
\nl
Let $\l_1,\l_2,\do,\l_k$ be a sequence elements of $\kk^*$ such that 

$\l_i-\l_j\ne0$ for $i\ne j$,

$\l_i-\l_j\i\ne0$ for all $i,j$.
\nl
For any $t\in[1,\s]$ the system of linear equations 
$$\sum_{i\in[1,p_t]}(\l_i^j-\l_i^{-j})c_{t,i}=-\d_{j,p_t}$$
($j\in[1,p_t]$) with unknowns $c_{t,i}$ ($i\in[1,p_t]$) has a unique solution 
$(c_{t,i})_{i\in[1,p_t]}\in\kk^{p_t}$. (Its determinant is nonzero by 1.1(a) with $m=p_t$.
Note that $\l_i$ is defined for $i\in[1,p_t]$ since $p_t\le p_1=k$.)
For any $i\in[1,k]$ we choose a basis $(v_{t,i})_{t\in[1,\s];p_t\ge i})$ of $\cv_i$ and we define vectors 
$v'_{t',i}\in\cv'_i$ ($t'\in[1,\s],p_{t'}\ge i$) by $(v_{t,i},v'_{t',i})=\d_{t,t'}c_{t,i}$ for all 
$t\in[1,\s],p_t\ge i$. Then for any $t\in[1,\s]$ and any $j\in[1,p_t]$ we have

(a) $\sum_{i\in[1,p_t]}(\l_i^j-\l_i^{-j})(v_{t,i},v'_{t,i})=-\d_{j,p_t}$.
\nl
Hence for $j\in[-p_t,p_t-1]$ we have

(b) $\sum_{i\in[1,p_t]}(\l_i^j-\l_i^{-j})(v_{t,i},v'_{t,i})=\d_{j,-p_t}$.
\nl
(For $j\in[1,p_t-1]$, (b) follows from (a); for $j\in[-p_t,-1]$, (b) follows from (a) by replacing $j$ by $-j$; 
for $j=0$, (b) is obvious.) 
For any $t\in[1,\s]$ we set $v_t=\sum_{i\in[1,k];i\le p_t}(v_{t,i}+v'_{t,i})\in V$. Define a linear map $g:V@>>>V$
by $gx=\l_ix$ for $x\in\cv_i$, $gx=\l_i\i x$ for $x\in\cv'_i$ ($i\in[1,k]$). Then $g\in G$.

\subhead 1.4\endsubhead
Assume that $Q\ne0$ so that $p\ne2$. Let $p_*=(p_1\ge p_2\ge\do\ge p_\s)$ be a sequence of integers $\ge1$ such 
that $p_1+p_2+\do+p_\s=n$. If $\k=0$ we assume also that $\k_\s=0$. 
For any $i\ge1$ we set $\bp_i=\sha(t\in[1,\s];p_t\ge i)$ so that $\bp_1\ge\bp_2\ge\do$ and $\sum_i\bp_i=n$. Let
$k=p_1$. We have $\bp_1=\s$, $\bp_k\ge1$, $\bp_{k+1}=0$. 

We can choose subspaces $\cz',\cz''$ of $V$ and subspaces $\cv_i,\cv'_i$ of $V$ (for $i\in[2,k]$) such that:

$V=\cz'\op\cz''\op\op_{i\in[2,k]}(\cv_i\op\cv'_i)$;

$(,)$ is nondegenerate on $\cz'$ and on $\cz''$, $()$ is zero on $\cv_i,\cv'_i$ (for $i\in[2,k]$);

$\dim\cz'=\bp_1+\k-\k_\s,\dim\cz''=\bp_1+\k_\s$;

$\dim\cv_i=\dim\cv'_i=\bp_i$ for $i\in[2,k]$;

$(\cv_i\op\cv'_i,\cv_j\op\cv'_j)=0$ for $i\ne j$ in $[2,k]$;

$(\cz',\cz'')=0$; 

$(\cz'+\cz'',\cv_i+\cv'_i)=0$ for $i\in[2,k]$.
\nl
Let $\l_i (i\in[2,k])$ be elements of $\kk^*$ such that 

$\l_i-\l_j\ne0$ for $i\ne j$ in $[2,k]$,

$\l_i-\l_j\i\ne0$ for all $i,j$ in $[2,k]$.
\nl
For any $t\in[1,\s]$, the system of linear equations 
$$c_{t,1}+(-1)^jc_{t,-1}+\sum_{i\in[2,p_t]}(\l_i^j+\l_i^{-j})c_{t,\l_i}=\d_{j,p_t}$$
($j\in[0,p_t]$) with $(p_t+1)$ unknowns ($c_{t,\l_i}\in\kk$ ($i\in[2,p_t]$) and $c_{t,1}\in\kk,c_{t,-1}\in\kk$) 
has a unique solution. (Its determinant is nonzero by 1.1(b) with $m=p_t-1$.)
For any $i\in[2,k]$ we choose a basis $(v_{t,i})_{t\in[1,\s];p_t\ge i})$ of $\cv_i$ and we define vectors 
$v'_{t',i}\in\cv'_i$ ($t'\in[1,\s],p_{t'}\ge i$) by $(v_{t,i},v'_{t',i})=\d_{t,t'}c_{t,\l_i}$ for all 
$t\in[1,\s],p_t\ge i$. (Note that if $t\in[1,\s]$ is such that $p_t\ge i$ then $i\in[2,p_t]$ hence
$c_{t,\l_i}$ is defined.) 

We can find vectors $v'_t\in\cz'$ ($t\in[1,\s]$) such that $(v'_t,v'_{t'})=c_{t,1}\d_{t,t'}$ for all $t,t'$ in 
$[1,\s]$. We can find vectors $v''_t\in\cz''$ ($t\in[1,\s]$) such that $(v''_t,v''_{t'})=c_{t,-1}\d_{t,t'}$ 
for all $t,t'$ in $[1,\s]$. Then for any $t\in[1,\s]$ and any $j\in[0,p_t]$ we have
$$(v'_t,v'_t)+(-1)^j(v''_t+v''_t)+\sum_{i\in[2,p_t]}(\l_i^j+\l_i^{-j})(v_{t,i},v'_{t,i})=\d_{j,p_t}.$$
It follows that
$$(v'_t,v''_t)+(-1)^j(v''_t,v''_t)+\sum_{i\in[2,p_t]}(\l_i^j+\l_i^{-j})(v_{t,i},v'_{t,i})=\d_{-j,p_t}$$
for $j\in[-p_t,p_t-1]$.
For any $t\in[1,\s]$ we set $v_t=v'_t+v''_t+\sum_{i\in[2,k];i\le p_t}(v_{t,i}+v'_{t,i})\in V$. Define a linear 
map $g:V@>>>V$ by $gx=x$ for $x\in\cz'$, $gx=-x$ for $x\in\cz''$, $gx=\l_ix$ for $x\in\cv_i$ ($i\in[2,k]$), 
$gx=\l_i\i x$ for $x\in\cv'_i$ ($i\in[2,k]$). Then $g\in G$.

\subhead 1.5\endsubhead
Assume that we are in the setup of 1.3 or 1.4. For $t,t'$ in $[1,\s]$ and $j\in[-p_t,p_t-1]$ we have
$$(g^jv_t,v_{t'})=\d_{t,t'}\d_{-j,p_t}.\tag a$$
Indeed, in the setup of 1.3, the left hand side of (a) is equal to
$$\align&\sum_{i\in[1,k];i\le p_t,i\le p_{t'}}(\l_i^jv_{t,i}+\l_i^{-j}v'_{t,i},v_{t',i}+v'_{t',i})
\\&=\d_{t,t'}\sum_{i\in[1,p_t]}(\l_i^j-\l_i^{-j})(v_{t,i},v'_{t,i})=\d_{t,t'}\d_{j,-p_t};\endalign$$
in the setup of 1.4, the left hand side of (a) is equal to
$$\align&(v_{t,1},v_{t',1})+(-1)^j(v'_{t,1},v'_{t',1})+
\sum_{i\in[2,k];i\le p_t,i\le p_{t'}}(\l_i^jv_{t,i}+\l_i^{-j}v'_{t,i},v_{t',i}+v'_{t',i})\\&
=\d_{t,t'}((v_{t,1},v_{t,1})+(-1)^j(v'_{t,1},v'_{t,1})\\&+
\sum_{i\in[2,k];i\le p_t,i\le p_t}(\l_i^j+\l_i^{-j})(v_{t,i},v'_{t,i}))=\d_{t,t'}\d_{j,-p_t}.\endalign$$
As in \cite{\WEUN, 3.3(vi)}, from (a) we deduce that the vectors $(g^jv_t)_{t\in[1,\s],j\in[-p_t,p_t-1]}$ span a 
$(\nn-\k)$-dimensional subspace of $V$ on which $(,)$ is nondegenerate. (If $\k=0$ this subspace is $V$.) For any
$h\in[1,n]$ we can write $h=p_1+p_2+\do+p_{r-1}+i$ where $r\in[1,\s]$ and $i\in[1,p_r]$ are uniquely determined; 
we define $V_h$ to be the subspace of $V$ spanned by the vectors $g^jv_t (t\in[1,r-1],j\in[0,p_t-1])$ and 
$g^jv_r (j\in[0,i-1])$. Let $V'_h$ be the subspace of $V$ spanned by the vectors 
$g^jv_t (t\in[1,r-1],j\in[1,p_t])$ and $g^jv_r (j\in[1,i])$. We have $(V_h,V_h)=0$ (see (a)), $(V'_h,V'_h)=0$, 
$gV_h=V'_h$. There are unique sequences $V_*,V'_*$ in $\cf$ such that $V_h,V'_h$ are as above for any $h\in[1,n]$.
We have $V'_*=gV_*$. For any $r\in[1,\s]$ and $i\in[1,p_r-1]$ we have
$$\dim(V'_{p_1+p_2+\do+p_{r-1}+i}\cap V_{p_1+p_2+\do+p_{r-1}+i})=p_1+p_2+\do+p_{r-1}+i-r,$$
(the intersection is spanned by the vectors $g^jv_t (t\in[1,r-1],j\in[1,p_t-1])$ and $g^jv_r (j\in[1,i-1])$),
$$\dim(V'_{p_1+p_2+\do+p_{r-1}+i}\cap V_{p_1+p_2+\do+p_{r-1}+i+1})=p_1+p_2+\do+p_{r-1}+i-r+1,$$
(the intersection is spanned by the vectors $g^jv_t (t\in[1,r-1],j\in[1,p_t-1])$ and $g^jv_r (j\in[1,i])$).
For any $r\in[1,\s]$ we have
$$\dim(V'_{p_1+p_2+\do+p_r}\cap V_{\nn-p_1-p_2-\do-p_{r-1}-1})=p_1+p_2+\do+p_r-r,$$
(the intersection is spanned by the vectors $g^jv_t (t\in[1,r],j\in[1,p_t-1])$),
$$\dim(V'_{p_1+p_2+\do+p_r}\cap V_{\nn-p_1-p_2-\do-p_{r-1}})=p_1+p_2+\do+p_r-r+1,$$
(the intersection is spanned by the vectors $g^jv_t (t\in[1,r],j\in[1,p_t-1])$ and $g^{p_r}v_r$). (We use again 
(a).) As in \cite{\WEUN, 3.2} we deduce that $a_{V_*,V'_*}=w_{p_*}\in\WW$ (notation of \cite{\WEUN, 1.4, 1.6}).
Let $B,B'$ be the stabilizers of $V_*,V'_*$ in $G$. Then $B,B'$ are Borel subgroups of $G$ and
$(B,B')\in\co_{w_{p_*}}$, $gBg\i=B'$. Hence if $\g$ denotes the conjugacy class of $g$ in $G$, we have 
$$\fB_{w_{p_*}}^\g\ne\em.\tag b$$ 
Note that $\g$ is a semisimple conjugacy class and that $w_{p_*}$ has minimal length in its conjugacy class $C$ in
$\WW$ (which is elliptic).

Let $\d(g)=\dim Z(g)$. Let $d=\ul(w_{p_*})$. We show that
$$\d(g)=d.\tag c$$
In the setup of 1.3, $Z(g)$ is isomorphic to $GL(\bp_1)\T GL(\bp_2)\T\do\T GL(\bp_k)$ hence 
$\d(g)=\bp_1^2+\bp_2^2+\do+\bp_k^2$. In the setup of 1.4, the identity component of $Z(g)$ is isomorphic to
$SO(\bp_1+\k-\k_\s)\T SO(\bp_1+\k_\s)\T GL(\bp_2)\T GL(\bp_3)\T\do\T GL(\bp_k)$ hence
$$\align&\d(g)=(\bp_1+\k-\k_s)(\bp_1+\k-\k_\s-1)/2+(\bp_1+\k_s)(\bp_1+\k_s-1)/2\\&+\bp_2^2+\do+\bp_k^2
=\bp_1^2+\bp_2^2+\do+\bp_k^2-\s(1-\k).\endalign$$ 
If $(1-\k)Q=0$ we have $d=2(p_2+2p_3+\do+(\s-1)p_\s)+n$; if $(1-\k)Q\ne0$ we have 
$d=2(p_2+2p_3+\do+(\s-1)p_\s)+n-\s$. Hence to prove (c) it is enough to show that
$$\bp_1^2+\bp_2^2+\do+\bp_k^2=2(p_2+2p_3+\do+(\s-1)p_\s)+n.$$
This follows from the equality $X=2Y$ in \cite{\WEUN, 4.4}; note that $f_{2h}$ from {\it loc.cit.} is the same as 
$\bp'$ and $\sum_hf_{2h}=n$. 
From (b) and (c) we see that $\g$ is a (semisimple) $C$-small conjugacy class. This proves 0.6(i) for our $G$.

\head 2. Exceptional groups\endhead
\subhead 2.1\endsubhead
In this subsection we assume that $\kk$ is an algebraic closure of a finite field $\FF_q$ with $q$ elements; we
also assume that 0.1(i) holds. We choose an $\FF_q$-split rational structure on $G$ with Frobenius map $F:G@>>>G$
such that $B^*$ and $\ct$ are $F$-stable. Note that $F(t)=t^q$ for all $t\in\ct$. Define a class function 
$\Pi_G:\WW@>>>\ZZ$ by $\Pi_G(w)=\sum_i\tr(w,H^{2i}(\cb,\bbq))q^i$ where we use the standard $\WW$-module 
structure on $H^*(\cb,\bbq)$. For any $z\in\cw$ let $x_z\in G$ be such that $x_z\i F(x_z)=z$ and let 
$\ct_z=x_z\ct x_z\i$, an $F$-stable maximal torus of $G$. For any $w\in\WW$ we have 
$$\Pi_G(w)=(-1)^{\ul(w)}|G^F|q^{-\nu_G}|\ct_w|\i.$$

\subhead 2.2\endsubhead
We assume that $\kk$ is as in 2.1, that $G$ is almost simple, simply connected of exceptional type and that 
$K\sneq\D$. Let $\g\in\cs_K$. Let $C\in\uWW_{el}$ and let $w\in C_{min}$. We show that the condition that 
$\fB^\g_w\ne\em$ can be tested by performing a computer calculation. We will also see that this condition depends
only on $K$, not on $\g$.

We choose an $\FF_q$-rational structure on $G$ as in 2.1. We can assume that $g^{q-1}=1$ for some/any $g\in\g$. 
Then $\g$ is $F$-stable and $\g\cap\ct=\g\cap\ct^F$ is a single $\cw$-orbit. Let $\z\in\g\cap\ct$ be such that 
$Z(\z)=G_K$. Note that $G_K$ is defined and split over $\FF_q$, 

Now the class function $\Pi_{G_K}:\WW_K@>>>\ZZ$ is well defined, see 2.1. For $z\in\cw$ we have
$$\align&\sha(h\in G^F;h\i\z h\in\ct_z)=\sha(h\in G^F;x_z\i h\i\z hx_z\in\ct)\\&=
\sha(h'\in G;F(h')=h'\dz,h'{}\i\z h'\in\ct)=\sha(h'\in G;F(h')=h'\dz,h'{}\i\z h'\in\g\cap\ct)\\&
=|\cw_K|\i\sum_{v\in\cw}\sha(h'\in G;F(h')=h'\dz,h'{}\i\z h'=\dv\z\dv\i)\\&
=|\cw_K|\i\sum_{v\in\cw}\sha(h''\in G;F(h'')=h''\dv\i\dz F(\dv),h''{}\i\z h''=\z)\\&
=|\cw_K|\i\sum_{v\in\cw}\sha(h''\in G_K;F(h'')=h''\dv\i\dz F(\dv))\\&
=|\cw_K|\i\sha(v\in\cw;\dv\i\dz F(\dv)\in G_K)|G_K^F|=|\cw_K|\i\sha(v\in\cw;v\i zv\in\cw_K)|G_K^F|.\endalign$$
(We set $h'=hx_z$; then we set $h'{}\i\z h'=\dv\z\dv\i$ with $v\in\cw$; then we set $h''=h'\dv$ and we use Lang's 
theorem in $G_K$.)

As in \cite{\WEUN, 1.2(a)} the number of fixed points of $F:\fB^\g_w@>>>\fB^\g_w$, $(g,B)\m(F(g),F(B))$, is given
by
$$|(\fB^\g_w)^F|=|\WW|\i\sum_{E,E'\in\Irr\WW,z\in\WW,g\in\g^F}\tr(T_w,E_q)(\r_E:R_{E'})\tr(z,E')\tr(\z,R^1(z)).
\tag a$$
(Notation of {\it loc.cit.}.) Using \cite{\DL, 7.2} we see that
$$\align&\tr(\z,R^1(z))=\sha(h\in G^F;h\i\z h\in\ct_z)|\ct_z|\i q^{-\nu_{G_K}}(-1)^{\ul(z)}\\&
=|\cw_K|\i\sha(v\in\cw;v\i zv\in\cw_K)|G_K^F|\ct_z|\i q^{-\nu_{G_K}}(-1)^{\ul(z)}\\&
=|\cw_K|\i\sha(v\in\cw;v\i zv\in\cw_K)\Pi_{G_K}(v\i zv).\endalign$$
(We have used that the restriction to $\cw_K=\WW_K$ of the function $z\m(-1)^{\ul(z)}$ on $\cw=\WW$ is the 
analogous function defined in terms of $G_K$.) Substituting this into (a) we obtain
$$\align&|(\fB^\g_w)^F|=|\g^F||\WW|\i\sum_{E,E'\in\Irr\WW,z\in\WW}\tr(T_w,E_q)(\r_E:R_{E'})\tr(z,E')\\&
\T|\cw_K|\i\sha(v\in\cw;v\i zv\in\cw_K)\Pi_{G_K}(v\i zv)\\&=|\g^F||\WW|\i
\sum_{E,E'\in\Irr\WW,z\in\WW}\tr(T_w,E_q)(\r_E:R_{E'})\tr(z,E')\tr(z,\ind_{\cw_K}^{\cw}(\Pi_{G_K})).\endalign$$ 
Hence 
$$|(\fB^\g_w)^F|=|G^F|/|G_K^F|\sum_{E,E'\in\Irr\WW}\tr(T_w,E_q)(\r_E:R_{E'})(E':\Pi_{G_K})_{\WW_K}.\tag b$$
Here $(E':\Pi_{G_K})_{\WW_K}$ is the inner product of $\Pi_{G_K}$ (viewed as a representation of $\WW_K$) with 
the restriction of $E'$ to $\WW_K$. We can also write (b) as follows:
$$|(\fB^\g_w)^F|=|G^F|/|G_K^F|\sum A_{E,C}\ph_{E,E'}m_{E',E''}t_{E'',K}$$
where the sum is taken over all $E,E'$ in $\Irr\WW$, $E''\in\Irr\WW_K$ and the notation is as follows. For 
$C'\in\uWW_{el},E\in\Irr\WW$ we set $A_{E,C'}=\tr(T_z,E_q)$ where $z\in C'_{min}$. (Note that $A_{E,C'}$ is well 
defined by \cite{\GP, 8.2.6(b)}.) For $E,E'\in\Irr\WW$ let $\ph_{E,E'}=(\r_E:R_{E'})$. (Notation of 
\cite{\WEUN, 1.2}.) For $E'\in\Irr\WW$, $E''\in\Irr\WW_K$ let $m_{E',E''}$ be the multiplicity of $E''$ in 
$E'|_{\WW_K}$. For $E''\in\Irr\WW_K$ let $t_{E'',K}$ be the multiplicity of $E''$ in $\Pi_{G_K}$. Thus 
$|(\fB^\g_w)^F|$ is $|G^F|/|G_K^F|$ times the the $C$-entry of the vector
$${}^t(A_{E,C'})(\ph_{E,E'})(m_{E',E''})(t_{E'',K}).$$
Here the matrix $(A_{E,C'})$ is known from the works of Geck and Geck-Michel (see \cite{\GP, 11.5.11}) and is 
available through the CHEVIE package \cite{\GH}. The matrix $\ph_{E,E'}$ has as entries the coefficients of the 
"nonabelian Fourier transform" in \cite{\OR, 4.15}. The matrix $(m_{E',E''})$ ("Induction table") and the vector 
$(t_{E'',K})$ ("Fake degree") are available through the CHEVIE package. Thus $|(\fB^\g_w)^F|$ can be obtained by 
calculating the product of several explicitly known matrices. The calculation was done  using the CHEVIE package.
It turns out that $|(\fB^\g_w)^F|$ is a polynomial in $q$ with integer coefficients denoted by $P^K_C$ (it 
depends only on $K,C$ not on $\g,w$). Note that $\fB^\g_w\ne\em$ if and only if $P^K_C\ne0$ as a polynomial in 
$q$. Thus the condition that $\fB^\g_w\ne\em$ can be tested. Moreover for each $K$ such that $P^K_C\ne0$ and for 
$\g\in\cs_K$, the condition
that $\g$ is $C$-small is equivalent to the condition that $\dim(G_K)=d_C$; in this case we have
$P^K_C=m^K_C|G(\FF_q)|$ as polynomials in $q$ where $m^K_C$ is an integer $\ge1$ independent of $q$. (For any $K$
such that $P^K_C\ne0$ we have $\deg(P^K_C)\ge\dim(G)$ (by \cite{\WEUN, 5.2}). Note that $m^K_C$ is equal to the 
number of connected components of $\fB^\g_w$ for $\g\in\cs_K$, $w\in C_{min}$. This number can be $>1$; in one 
example in type $E_8$ it is $10$.

\subhead 2.3\endsubhead
In this subsection we give (in the setup of 2.2) tables which describe for each exceptional type and each 
$C\in\uWW_{el}$ (with one exception) some proper subsets $K$ of $\D$ such that $P^K_C\ne0$ and $\dim(G_K)=d_C$. 
The elements of $\D-\{\a_0\}$ are denoted by numbers $1,2,3,\do$ as in \cite{\GP, p.20}. We write $0$ instead of 
$\a_0$. We specify $K$ by marking each element of $\D-K$ by $\bul$.
An element $C\in\uWW_{el}$ is specified by indicating the characteristic polynomial of an element of $C$ acting 
on $R_\WW$, a product of cyclotomic polynomials $\Ph_d$ (an exception is type $F_4$
when there are two choices for $C$ with characteristic polynomial $\Ph_2^2\Ph_6$ in which case we use the 
notation $(\Ph_2^2\Ph_6)'$, $(\Ph_2^2\Ph_6)''$ for what in \cite{\GP, p.407} is denoted by $D_4$, $C_3+A_1$). The
notation $d;C;\c;(K_1)_{m_1};(K_2)_{m_2};\do$ means that $C\in\uWW_{el}$, $d=d_C$, $\c=\r_{\g_C}$ ($\g_C$ as in 
0.4), and $K_1,K_2,\do$ are proper subsets of $\D$ such that $P^{K_i}_C\ne0$ and $\dim(G_{K_i})=d_C$; we have 
$m_i=m^{K_i}_C$. (We omit $m_i$ whenever $m_i=1$.)

The notation for irreducible representations of $\WW$ (of type $E_6,E_7,E_8$) is as in \cite{\SPS}; for type 
$F_4$ it is as in \cite{\OR}; for type $G_2$, $1_0$ is the unit representation, $2_1$ is the reflection 
representation and $2_2$ is the other two dimensional irreducible representation of $\WW$.

\mpb

{\it Type} $G_2$; $\D$ is $(012)$

$2;\Ph_6; 1_0;\qua (\bul\bul\bul)$
          
$4;\Ph_3; 2_1;\qua (\bul1\bul);\qua (\bul\bul2)$

$6;\Ph_2^2; 2_2;\qua (0\bul2)$

\mpb

{\it Type} $F_4$; $\D$ is $(01234)$

$4;\Ph_{12};  1_1 ; \qua (\bul\bul\bul\bul\bul)$

$6; \Ph_8;  4_2 ; \qua(\bul1\bul\bul\bul)$

$8; \Ph_6^2;  9_1; \qua(0\bul2\bul\bul);\qua (\bul1\bul3\bul)_2$

$10; (\Ph_2^2\Ph_6)'; 8_1 ; \qua(\bul\bul\bul 3 4)$

$10; (\Ph_2^2\Ph_6)'';  8_3 ; \qua(\bul12\bul\bul);\qua (0\bul2\bul4)$

$12;\Ph_4^2;12_1 ;\qua(\bul\bul23\bul)_3;\qua(\bul12\bul4);\qua(0\bul\bul34)$

$14; \Ph_2^2\Ph_4; 16_1 ;\qua(0\bul23\bul)$

$16;  \Ph_3^2; 6_1  ; \qua(01\bul34)$

$24;  \Ph_2^4; 9_4  ; \qua(0\bul234)$

\mpb

{\it Type} $E_6$; $\D$ is $\left(\sm 1&3&4&5&6\\{}&{}&2&{}&{}\\{}&{}&0&{}&{}\esm\right)$

$6;   \Ph_3\Ph_{12}; 1_0; \qua\left(\sm\bul&\bul&\bul&\bul&\bul\\{}&{}&\bul&{}&{}\\{}&{}&\bul&{}&{}\esm\right)$

$8;  \Ph_9; 6_1;\qua\left(\sm1&\bul&\bul&\bul&\bul\\{}&{}&\bul&{}&{}\\{}&{}&\bul&{}&{}\esm\right)$

$12; \Ph_3\Ph_6^2; 30_3; \qua\left(\sm1&3&\bul&\bul&\bul\\{}&{}&\bul&{}&{}\\{}&{}&\bul&{}&{}\esm\right);
\qua\left(\sm 1&\bul&\bul&5&\bul\\{}&{}&2&{}&{}\\{}&{}&\bul&{}&{}\esm\right)$

$14; \Ph_2^2\Ph_3\Ph_6; 15_4;\qua\left(\sm1&\bul&4&\bul&6\\{}&{}&\bul&{}&{}\\{}&{}&0&{}&{}\esm\right)$

$24; \Ph_3^3; 10_9;\qua\left(\sm1&3&\bul&5&6\\{}&{}&2&{}&{}\\{}&{}&0&{}&{}\esm\right)$

\mpb

{\it Type} $E_7$; $\D$ is $\left(\sm0&1&3&4&5&6&7\\{}&{}&{}&2&{}&{}&{}\esm\right)$

$7; \Ph_2\Ph_{18}; 1_0;\qua\left(\sm\bul&\bul&\bul&\bul&\bul&\bul&\bul\\{}&{}&{}&\bul&{}&{}&{}\esm\right)$

$9; \Ph_2\Ph_{14}; 7_1; \qua\left(\sm\bul&1&\bul&\bul&\bul&\bul&\bul\\{}&{}&{}&\bul&{}&{}&{}\esm\right)$

$11; \Ph_2\Ph_6\Ph_{12}; 27_2; \qua\left(\sm\bul&1&\bul&\bul&\bul&\bul&\bul\\{}&{}&{}&2&{}&{}&{}\esm\right)$

$13; \Ph_2\Ph_6\Ph_{10}; 56_3; \qua\left(\sm\bul&1&3&\bul&\bul&\bul&\bul\\{}&{}&{}&\bul&{}&{}&{}\esm\right)_2;
\qua\left(\sm\bul&1&\bul&\bul&5&\bul&\bul\\{}&{}&{}&2&{}&{}&{}\esm\right)$

$15;  \Ph_2^3\Ph_{10}; 35_4; \qua\left(\sm0&\bul&3&\bul&5&\bul&7\\{}&{}&{}&\bul&{}&{}&{}\esm\right)$

$17;  \Ph_2\Ph_4\Ph_8; 189_5; \qua\left(\sm\bul&1&3&\bul&5&\bul&\bul\\{}&{}&{}&2&{}&{}&{}\esm\right); \qua
\left(\sm0&\bul&3&\bul&5&\bul&7\\{}&{}&{}&2&{}&{}&{}\esm\right)$

$21; \Ph_2\Ph_6^3; 315_7; \qua\left(\sm\bul&1&\bul&4&5&\bul&\bul\\{}&{}&{}&2&{}&{}&{}\esm\right)_3;
\qua\left(\sm\bul&1&3&\bul&5&6&\bul\\{}&{}&{}&2&{}&{}&{}\esm\right)$

$23; \Ph_2^3\Ph_6^2; 280_8;\qua\left(\sm0&1&3&\bul&5&\bul&7\\{}&{}&{}&\bul&{}&{}&{}\esm\right)$

$25; \Ph_2\Ph_3^2\Ph_6; 70_9;\qua\left(\sm0&1&\bul&4&\bul&6&7\\{}&{}&{}&2&{}&{}&{}\esm\right)$

$31; \Ph_2^5\Ph_6; 84_{12}; \qua\left(\sm0&1&3&\bul&5&6&7\\{}&{}&{}&\bul&{}&{}&{}\esm\right)$

$33; \Ph_2^3\Ph_4^2; 210_{13}; \qua\left(\sm\bul&1&3&4&\bul&6&7\\{}&{}&{}&2&{}&{}&{}\esm\right)$

$63; \Ph_2^7; 15_{28}; \qua\left(\sm0&1&3&4&5&6&7\\{}&{}&{}&\bul&{}&{}&{}\esm\right)$

\mpb

{\it Type} $E_8$; $\D$ is $\left(\sm1&3&4&5&6&7&8&0\\{}&{}&2&{}&{}&{}&{}&{}\esm\right)$

$8; \Ph_{30}; 1_0; \qua\left(\sm\bul&\bul&\bul&\bul&\bul&\bul&\bul&\bul\\{}&{}&\bul&{}&{}&{}&{}&{}\esm\right)$

$10; \Ph_{24}; 8_1; \qua\left(\sm1&\bul&\bul&\bul&\bul&\bul&\bul&\bul\\{}&{}&\bul&{}&{}&{}&{}&{}\esm\right)$

$12; \Ph_{20}; 35_2; \qua\left(\sm1&\bul&\bul&\bul&\bul&\bul&\bul&\bul\\{}&{}&2&{}&{}&{}&{}&{}\esm\right)$

$14; \Ph_6\Ph_{18}; 112_3;\qua\left(\sm1&3&\bul&\bul&\bul&\bul&\bul&\bul\\{}&{}&\bul&{}&{}&{}&{}&{}\esm\right)_2;
\qua\left(\sm1&\bul&\bul&5&\bul&\bul&\bul&\bul\\{}&{}&2&{}&{}&{}&{}&{}\esm\right)$

$16; \Ph_{15}; 210_4; \qua\left(\sm1&\bul&4&\bul&\bul&\bul&\bul&\bul\\{}&{}&2&{}&{}&{}&{}&{}\esm\right)_2;
\left(\sm1&\bul&\bul&5&\bul&7&\bul&\bul\\{}&{}&2&{}&{}&{}&{}&{}\esm\right)$

$16; \Ph_2^2\Ph_{18}; 84_4; \qua\left(\sm\bul&\bul&\bul&5&\bul&7&\bul&0\\{}&{}&2&{}&{}&{}&{}&{}\esm\right)$

$18; \Ph_2^2\Ph_{14} ; 560_5;\qua\left(\sm1&3&\bul&5&\bul&\bul&\bul&\bul\\{}&{}&2&{}&{}&{}&{}&{}\esm\right);
\qua\left(\sm1&\bul&\bul&5&\bul&7&\bul&0\\{}&{}&2&{}&{}&{}&{}&{}\esm\right)$

$20; \Ph_{12}^2 ; 700_6; \qua\left(\sm1&3&\bul&5&6&\bul&\bul&\bul\\{}&{}&\bul&{}&{}&{}&{}&{}\esm\right)_2;
\qua\left(\sm1&3&\bul&5&\bul&7&\bul&\bul\\{}&{}&2&{}&{}&{}&{}&{}\esm\right)$

$22; \Ph_4^2\Ph_{12}; 400_7; \qua\left(\sm1&3&\bul&5&\bul&7&\bul&0\\{}&{}&2&{}&{}&{}&{}&{}\esm\right)$

$22; \Ph_6^2\Ph_{12}; 1400_7;\qua\left(\sm1&\bul&4&5&\bul&\bul&\bul&\bul\\ {}&{}&2&{}&{}&{}&{}&{}\esm\right)_3;
\qua\left(\sm1&3&\bul&5&6&\bul&\bul&\bul\\{}&{}&2&{}&{}&{}&{}&{}\esm\right)$

$24; \Ph_{10}^2; 1400_8;\qua\left(\sm1&\bul&4&5&\bul&7&\bul&\bul\\{}&{}&2&{}&{}&{}&{}&{}\esm\right)_3;
\qua\left(\sm1&3&\bul&5&6&\bul&8&\bul\\{}&{}&2&{}&{}&{}&{}&{}\esm\right)$

$24; \Ph_2^2\Ph_6\Ph_{12}; 1344_8;\qua\left(\sm\bul&\bul&4&5&\bul&7&\bul&0\\ {}&{}&2&{}&{}&{}&{}&{}\esm\right)$

$26; \Ph_3^2\Ph_{12}; 448_9;\qua\left(\sm1&3&\bul&5&6&\bul&8&0\\{}&{}&\bul&{}&{}&{}&{}&{}\esm\right)$

$26;  \Ph_2^2\Ph_6\Ph_{10}; 3240_9; \qua\left(\sm1&3&4&\bul&6&7&\bul&\bul\\{}&{}&\bul&{}&{}&{}&{}&{}\esm\right);
\qua\left(\sm1&\bul&4&5&\bul&7&\bul&0\\{}&{}&2&{}&{}&{}&{}&{}\esm\right)$

$28; \Ph_3\Ph_9; 2240_{10}; \qua\left(\sm1&3&\bul&5&6&7&\bul&\bul\\{}&{}&2&{}&{}&{}&{}&{}\esm\right);
\qua\left(\sm1&3&\bul&5&6&\bul&8&0\\{}&{}&2&{}&{}&{}&{}&{}\esm\right)$

$30; \Ph_8^2; 1400_{11}; \qua\left(\sm1&3&\bul&5&6&7&\bul&0\\{}&{}&2&{}&{}&{}&{}&{}\esm\right)$

$32; \Ph_2^4\Ph_{10}; 972_{12}; \qua\left(\sm\bul&\bul&4&5&\bul&7&8&0\\{}&{}&2&{}&{}&{}&{}&{}\esm\right)$

$34; \Ph_2^2\Ph_4\Ph_8; 4536_{13}; \qua\left(\sm1&3&4&\bul&6&7&\bul&\bul\\ {}&{}&2&{}&{}&{}&{}&{}\esm\right);
\qua\left(\sm1&\bul&4&5&\bul&7&8&0\\ {}&{}&2&{}&{}&{}&{}&{}\esm\right)$

$40; \Ph_6^4; 4480_{16}; \qua\left(\sm1&\bul&4&5&6&7&\bul&\bul\\{}&{}&2&{}&{}&{}&{}&{}\esm\right)_{10};
\qua\left(\sm1&3&4&5&\bul&7&8&0\\ {}&{}&\bul&{}&{}&{}&{}&{}\esm\right)$

$42;  \Ph_2^2\Ph_6^3; 7168_{17}; \qua\left(\sm1&\bul&4&5&6&7&\bul&0\\{}&{}&2&{}&{}&{}&{}&{}\esm\right)$

$44; \Ph_2^4\Ph_6^2; 4200_{18}; \qua\left(\sm\bul&3&4&5&\bul&7&8&0\\{}&{}&2&{}&{}&{}&{}&{}\esm\right)$

$44; \Ph_3^2\Ph_6^2; 3150_{18}; \qua\left(\sm1&3&4&5&6&\bul&8&0\\{}&{}&\bul&{}&{}&{}&{}&{}\esm\right)$

$46; \Ph_2^2\Ph_3^2\Ph_6; 2016_{19}; \qua\left(\sm1&3&\bul&5&6&7&8&0\\{}&{}&2&{}&{}&{}&{}&{}\esm\right)$

$46; \Ph_2^2\Ph_4^2\Ph_6; 1344_{19};$

$48; \Ph_5^2; 420_{20}; \qua\left(\sm1&3&4&\bul&6&7&8&0\\{}&{}&2&{}&{}&{}&{}&{}\esm\right)$

$60; \Ph_4^4; 840_{26}; \qua\left(\sm1&3&4&5&\bul&7&8&0\\{}&{}&2&{}&{}&{}&{}&{}\esm\right)$

$64; \Ph_2^6\Ph_6; 700_{28};\qua\left(\sm\bul&\bul&4&5&6&7&8&0\\{}&{}&2&{}&{}&{}&{}&{}\esm\right)$

$66; \Ph_2^4\Ph_4^2; 1400_{29};\qua\left(\sm1&\bul&4&5&6&7&8&0\\{}&{}&2&{}&{}&{}&{}&{}\esm\right)$

$80; \Ph_3^4; 175_{36}; \qua\left(\sm1&3&4&5&6&7&8&0\\{}&{}&\bul&{}&{}&{}&{}&{}\esm\right)$

$120; \Ph_2^8; 50_{56}; \qua\left(\sm\bul&3&4&5&6&7&8&0\\{}&{}&2&{}&{}&{}&{}&{}\esm\right)$

\subhead 2.4\endsubhead
We prove 0.6(i) in the case where $G$ is almost simple, simply connected of exceptional type. Let 
$C\in\uWW_{el}$. Let $K$ be a proper subset of $\D$ associated to $C$ in the tables in 2.3. We can find
$\g\in\cs_K$ such that 

(i) any element of $\g$ has finite order.
\nl
We show that $\g$ is $C$-small. By a standard argument we are reduced to the case where $\kk$ is as in 2.1. In 
this case the calculations outlined in 2.2 show that $\g$ is $C$-small, as claimed. 

Next we consider for $w\in C_{min}$, the map $\p:\fB_w@>>>G$, $(g,B)\m G$. Let $\fK=\p_!\bbq$. Using 
\cite{\CSIII, 14.2(a)}, we see that the cohomology sheaves of $\fK$ behave smoothly when restricted to
$\cup_{\g\in\cs_K}\g$ (which is one of the pieces $Y_{L,\S}$ in \cite{\CSI, 13.11}). Since we know that when 
$\g\in\cs_K$ satisfies (i), some cohomology sheaf of $\fK$ is non-zero on $\g$, it follows that for any
$\g\in\cs_K$, some cohomology sheaf of $\fK$ is non-zero on $\g$; in particular, $\fB_w^\g\ne\em$. It follows 
that any $\g\in\cs_K$ is $C$-good. This completes the proof of Theorem 0.6(i).

\head 3. Springer representations\endhead
\subhead 3.1\endsubhead
In the setup of 1.2 we assume that $G$ is the identity component of $Is(V)$. 
If $(1-\k)Q=0$ we identify, as in \cite{\WEUN, 1.5}, $\WW$ with $W$, the group of all permutations of $[1,\nn]$ 
that commute with the involution $i\m\nn+1-i$. Let $W_n$ be the group of all permutations of $[1,2n]$ that 
commute with the involution $i\m2n+1-i$. If $\nn=2n$ we have  $W=W_n$; if $\nn=2n+1$ we identify $W$ with $W_n$ 
(hence $\WW$ with $W_n$) by $w\m w'$ where $h(w'(i))=w(h(i))$ for $i\in[1,2n]$ and $h(i)=i$ if $i\in[1,n]$,
$h(i)=i+1$ if $i\in[n+1,2n]$. As in \cite{\OR, 4.5} we write the irreducible representations of $W_n$ in the form
$[(\l_1>\l_2>\do>\l_{m+1}),(\mu_1>\mu_2>\do>\mu_m)]$ where $\l_i,\mu_i\in\NN$, $\sum_i\l_i+\sum_i\mu_i=m^2+n$ and
$m$ is sufficiently large. For example $[(0<1<2<\do<n),(1<2<\do<n)]$ is the sign representation of $W_n$.

If $(1-\k)Q\ne0$ we identify as in \cite{\WEUN, 1.5} $\WW$ with $W'_n$, the group of even permutations in $W_n$.
As in \cite{\OR, 4.6} we write the irreducible representations of $W'_n$ as unordered pairs
$[(\l_1>\l_2>\do>\l_m),(\mu_1>\mu_2>\do>\mu_m)]$ where $\l_i,\mu_i\in\NN$, $\sum_i\l_i+\sum_i\mu_i=m^2-m+n$  and 
$m$ is sufficiently large. (There are two irreducible representations corresponding to 
$[(\l_1>\l_2>\do>\l_m),(\mu_1>\mu_2>\do>\mu_m)]$ with $\l_i=\mu_i$ for all $i$.) For example 
$[(1<2<\do<n),(0<1<2<\do<n-1)]$ is the sign representation of $W'_n$.

\subhead 3.2\endsubhead
Let $S_n$ be the symmetric group in $n$ letters. Using \cite{\SPE, 5.3, 4.4(a)} we see that

(a) if $n=2c\in2\NN$ then $W_c\T W'_c,S_n$ are naturally reflection subgroups of $W_n$ and we have
$$\align&j_{W_c\T W'_c}^{W_n}(\sg)=j_{S_n}^{W_n}(\sg)\\&=[(0<2<3<4<\do<c+1),(1<2<3<\do<c)];\endalign$$
(b) if $n=2c+1\in2\NN+1$ then $W_c\T W'_{c+1},S_n$ are naturally reflection subgroups of $W_n$ and we have 
$$\align&j_{W_c\T W'_{c+1}}^{W_n}(\sg)=j_{S_n}^{W_n}(\sg)\\&=[(1<2<3<4<\do<c+1),(1<2<3<\do<c)].\endalign$$
Using \cite{\SPE, 6.3} we see that

(c) if $n=2c\in2\NN$ then $W'_c\T W'_c$ is naturally a reflection subgroup of $W'_n$ and we have
$$j_{W'_c\T W'_c}^{W'_n}(\sg)=[(2<3<4<\do<c+1),(0<1<2<3<\do<c-1)].$$
Let $p_1\ge p_2\ge\do\ge p_\s$ be integers $\ge1$ such that $p_1+\do+p_\s=n$. Define $\bp_1\ge\bp_2\ge\do\ge\bp_k$
as in 1.3. Note that $\bp_1=\s,k=p_1$. Define $\tp_1\ge\tp_2\ge\do\ge\tp_\s$ by $\tp_i=p_i-1$ if $i\in[1,\bp_k]$, 
$\tp_i=p_i$ if $i\in[\bp_k+1,\s]$. 

Assuming that $k>1$ we have $\tp_\s=p_\s\ge1$ and using \cite{\SPE, 4.4(a)} we see that

(d) if $\s=2\t+1$ we have 
$$\align&[(p_\s<p_{\s-2}+1<\do<p_3+\t-1<p_1+\t),(p_{\s-1}<p_{\s-3}+1<\do<p_2+\t-1)]\\&=
j_{S_{\bp_k}\T W_{n-\bp_k}}^{W_n}
(\sg\bxt[(\tp_\s<\tp_{\s-2}+1<\do<\tp_3+\t-1<\tp_1+\t),\\&(\tp_{\s-1}<\tp_{\s-3}+1<\do<\tp_2+\t-1)]);\endalign$$
(e) if $\s=2\t$ we have 
$$\align&[(0<p_{\s-1}+1<p_{\s-3}+2<\do<p_3+\t-1<p_1+\t),\\&(p_\s<p_{\s-2}+1<\do<p_2+\t-1)]\\&=
j_{S_{\bp_k}\T W_{n-\bp_k}}^{W_n}(\sg\bxt[(0<\tp_{\s-1}+1<\tp_{\s-3}+2<\do<\tp_3+\t-1<\tp_1+\t),\\&
(\tp_\s<\tp_{\s-2}+1<\do<\tp_2+\t-1)]).\endalign$$
Assuming that $k>1$, $\s=2\t$ we have $\tp_\s=p_\s\ge1$, $n-\bp_k\ge2$ and using \cite{\SPE, 6.2(a)} we see that
$$\align&[(p_{\s-1}+1<p_{\s-3}+2<\do<p_3+\t-1<p_1+\t),\\&(p_\s-1<p_{\s-2}<\do<p_4+\t-3<p_2+\t-2)]\\&=
j_{S_{\bp_k}\T W'_{n-\bp_k}}^{W'_n}(\sg\bxt[(\tp_{\s-1}+1
<\tp_{\s-3}+2<\do<\tp_1+\t),\\&(\tp_\s-1<\tp_{\s-2}<\do<\tp_4+\t-3<\tp_2+\t-2)]).\tag f\endalign$$
We show that

(g) $j_{S_{\bp_k}\T\do\T S_{\bp_2}\T S_{\bp_1}}^{W_n}(\sg)$ is equal to
$$[(p_\s<p_{\s-2}+1<\do<p_3+\t-1<p_1+\t),(p_{\s-1}<p_{\s-3}+1<\do<p_2+\t-1)]$$
if $\s=2\t+1$ and to
$$\align&[(0<p_{\s-1}+1<p_{\s-3}+2<\do<p_3+\t-1<p_1+\t),\\&(p_\s<p_{\s-2}+1<\do<p_2+\t-1)]\endalign$$
if $\s=2\t$. We argue by induction on $k$. If $k=1$ we have $\s=n$ and the result follows from (a),(b). If $k>1$ 
then 
$$j_{S_{\bp_k}\T\do\T S_{\bp_2}\T S_{\bp_1}}^{W_n}(\sg)=j_{S_{\bp_k}\T W_{n-\bp_k}}^{W_n}
(\sg\bxt j_{S_{\bp_{k-1}}\T\do\T S_{\bp_2}\T S_{\bp_1}}^{W_{n-\bp_k}}(\sg))$$ 
and the result follows from (d),(e) using the induction hypothesis and the transitivity of the $j$-induction.

We write $\bp_1=a+b$ where $b-a\in\{0,1\}$. We show that

(h) $j_{S_{p'_k}\T\do\T S_{p'_2}\T W_a\T W'_b}^{W_n}(\sg)$ is equal to

$[(p_\s<p_{\s-2}+1<\do<p_3+\t-1<p_1+\t),(p_{\s-1}<p_{\s-3}+1<\do<p_2+\t-1)]$ if $\s=2\t+1$ and to
$$\align&[(0<p_{\s-1}+1<p_{\s-3}+2<\do<p_3+\t-1<p_1+\t),\\&(p_\s<p_{\s-2}+1<\do<p_2+\t-1)]\endalign$$ 
if $\s=2\t$. 
\nl
Using (a),(b) and the transitivity of $j$-induction we see that
$$\align&j_{S_{\bp_k}\T\do\T S_{\bp_2}\T W_a\T W'_b}^{W_n}(\sg)
=j_{S_{\bp_k}\T\do\T\do\T S_{\bp_2}\T W_{\bp_1}}^{W_n}
(j_{S_{\bp_k}\T\do\T S_{\bp_2}\T W_a\T W'_b}^{S_{\bp_k}\T\do\T S_{\bp_2}\T W_{\bp_1}}(\sg))\\&
=j_{S_{\bp_k}\T\do\T S_{\bp_2}\T W_{\bp_1}}^{W_n}
(j_{S_{\bp_k}\T\do \T S_{\bp_2}\T S_{\bp_1}}^{S_{\bp_k}\T\do\T S_{\bp_2}\T W_{\bp_1}}(\sg))
=j_{S_{\bp_k}\T\do\T S_{\bp_2}\T S_{\bp_1}}^{W_n}(\sg))\endalign$$ 
and it remains to use (g).

Assuming that $\s=2\t$ we show that
$$\align&j_{S_{\bp_k}\T\do\T S_{\bp_2}\T W'_{\bp_1/2}\T W'_{\bp_1/2}}^{W'_n}(\sg)=
[(p_{\s-1}+1<p_{\s-3}+2<\do<p_1+\t),\\&(p_\s-1<p_{\s-2}<\do<p_4+\t-3<p_2+\t-2)].\tag i\endalign$$
We argue by induction on $k$. If $k=1$ we have $\s=n$ and the result follows from (c). If $k>1$ then the left 
hand side of (i) is equal to 
$$j_{S_{\bp_k}\T W'_{n-\bp_k}}^{W'_n}
(\sg\bxt j_{S_{\bp_{k-1}}\T\do\T S_{\bp_2}\T W'_{\bp_1/2}\T W'_{\bp_1/2}}^{W'_{n-\bp_k}}(\sg))$$
and the result follows from (f) using the induction hypothesis and the transitivity of the $j$-induction.

\subhead 3.3\endsubhead
Assume that we are in the setup of 1.3. Let $p_*=(p_1\ge p_2\ge\do\ge p_\s)$ be as in 1.3. We consider a 
unipotent class $\g$ in $G$ such that any $u\in\g$ has Jordan blocks of sizes 

(i) $2p_1,2p_2,\do,2p_\s$.
\nl
We set $\s=2\t+\k_\s$. We show:
$$\align&\r_\g=[(0<p_{\s-1}+1<p_{\s-3}+2<\do<p_3+\t-1<p_1+\t),\\&(p_\s<p_{\s-2}+1<\do<p_4+\t-2<p_2+\t-1)]
\text{ if $\k_\s=0$},\\&
\r_\g=[(p_\s<p_{s-2}+1<\do<p_3+\t-1<p_1+\t),\\&(p_{\s-1}<p_{\s-3}+1<\do< p_4+\t-2<p_2+\t-1)]
\text{ if $\k_\s=1$}.\tag a\endalign$$
To the partition (i) we will apply the procedure of \cite{\IC, 11.6}. Let $M=\s+\k_\s$. Let 
$z_M\ge\do\ge z_2\ge z_1$ be the sequence (i) if $\k_\s=0$ and $2p_1,2p_2,\do,2p_\s,0$ if $\k_\s=1$.  The 
sequence $z'_M>\do>z'_2>z'_1$ in {\it loc.cit.} is 

$2p_1+\s-1,2p_2+\s-2,\do,2p_\s$ (if $\k_\s=0$),

$2p_1+\s,2p_2+\s-1,\do,2p_\s+1,0$ (if $\k_\s=1$).
\nl
This contains $M/2$ even numbers $2y_{M/2}>\do>2y_2>2y_1$ given by

$\{2p_t+\s-t;t\in[1,\s],\k_t=0\}$ (if $\k_\s=0$), 

$\{2p_t+\s-t+1;t\in[1,\s],\k_t=0\}\sqc\{0\}$ (if $\k_\s=1$)
\nl
and $M/2$ odd numbers $2y'_{M/2}+1>\do>2y'_2+1>2y'_1+1$ given by

$\{2p_t+\s-t+\k_\s;t\in[1,\s],\k_t=1\}$.
\nl
Thus, the sets $(\{y'_{M/2}>\do>y'_2>y'_1\},\{y_{M/2}>\do>y_2>y_1\})$ are given by

$(\{p_t+\t-(t+1)/2;t\in[1,\s],\k_t=1\},\{p_t+\t-t/2;t\in[1,\s],\k_t=0\})$ (if $\k_\s=0$)

$(\{p_t+\t+(1-t)/2;t\in[1,\s],\k_t=1\}, \{p_t+\t-t/2+1;t\in[1,\s],\k_t=0\}\sqc\{0\})$ (if $\k_\s=1$).
\nl
If $\k_\s=0$, the multisets 

$(\{y'_\t-(\t-1)\ge\do\ge y'_2-1\ge y'_1\ge 0\},\{y_\t-(\t-1)>\do>y_2-1>y_1\})$ 
\nl
are given by

$(\{p_1\ge p_3\ge\do\ge p_{\s-1}\ge 0\},\{p_2\ge p_4\ge\do\ge p_\s\})$.
\nl
If $\k_\s=1$, the multisets 

$(\{y'_{\t+1}-\t\ge\do\ge y'_2-1\ge y'_1\ge 0\},\{y_{\t+1}-\t\ge\do\ge y_2-1\ge y_1\})$
\nl
are given by

$(\{p_1\ge p_3\ge\do\ge p_\s\ge 0\},\{p_2\ge p_4\ge\do\ge p_{\s-1}\ge0\})$.
\nl
Now (a) follows from \cite{\IC, \S12}. Using (a) and 3.2(g) we see that
$$\r_g=j_{S_{\bp_k}\T\do\T S_{\bp_2}\T S_{\bp_1}}^{W_n}(\sg).\tag b$$

\subhead 3.4\endsubhead
Assume that we are in the setup of 1.4. Let $p_*=(p_1\ge p_2\ge\do\ge p_\s)$ be as in 1.4. Define 
$\ps:[1,\s]@>>>\{-1,0,1\}$ by $\ps(t)=1$ if $t$ is odd and $p_{t-1}>p_t$ (the last condition is regarded as 
satisfied when $t=1$); $\ps(t)=-1$ if $t$ is even and $p_t>p_{t+1}$ (the last condition is regarded as satisfied 
when $t=\s$); $\ps(t)=0$ for all other $t$. We set $\s=2\t+\k_\s$. If $\nn=2n$ we assume that $\s=2\t$. We 
consider a unipotent class $\g$ in $G$ such that any $u\in\g$ has Jordan blocks of sizes

(i) $2p_1+\ps(1),2p_2+\ps(2),\do,2p_\s+\ps(\s)$ if $\nn=2n$ (hence $\k_\s=0$),

(ii) $2p_1+\ps(1),2p_2+\ps(2),\do,2p_\s+\ps(\s)$ if $\nn=2n+1$ and $\k_\s=1$,

(iii) $2p_1+\ps(1),2p_2+\ps(2),\do,2p_\s+\ps(\s),1$ if $\nn=2n+1$ and $\k_\s=0$.
\nl
We show:
$$\align&\r_\g=[(p_\s-1<p_{\s-2}<\do<p_4+\t-3<p_2+\t-2),\\&(p_{\s-1}+1<p_{\s-3}+2<\do<p_3+\t-1<p_1+\t)],
\text{ in case (i)},\\&
\r_\g=[(p_\s<p_{\s-2}+1<\do<p_3+\t-1<p_1+\t),\\&(p_{\s-1}<p_{\s-3}+1<\do< p_4+\t-2<p_2+\t-1)]
\text{ in case (ii)},\\&\r_\g=[0<p_{\s-1}+1<p_{\s-3}+2<\do<p_3+\t-1<p_1+\t),
\\&(p_\s<p_{\s-2}+1<\do<p_4+\t-2<p_2+\t-1)]\text{ in case (iii)}.\tag a\endalign$$
To the partition (i),(ii) or (iii) we will apply the procedure of \cite{\IC, 11.7}. Let $z_M\ge\do\ge z_2\ge z_1$
be the sequence (i),(ii) or (iii) (where $M=\s$ in cases (i),(ii) and $M=\s+1$ in case (iii)). The sequence 
$z'_M>\do>z'_2>z'_1$ in {\it loc.cit.} is

$2p_1+\ps(1)+\s-1,2p_2+\ps(2)+\s-2,\do,2p_\s+\ps(\s)$ (in cases (i),(ii)), 

$2p_1+\ps(1)+\s,2p_2+\ps(2)+\s-1,\do,2p_\s+\ps(\s)+1,1$ (in case (iii)).
\nl
This contains $[M/2]$ even numbers $2y_{[M/2]}>\do>2y_2>2y_1$ given by

$\{2p_t+\ps(t)+\s-t;t\in[1,\s],\k_t=\k_{\ps(t)}\}$ in case (i),

$\{2p_t+\ps(t)+\s-t;t\in[1,\s],\k_t\ne\k_{\ps(t)}\}$ in case (ii),

$\{2p_t+\ps(t)+\s-t+1;t\in[1,\s],\k_t\ne k_{\ps(t)}\}$ in case (iii),
\nl
and $[(M+1)/2]$ odd numbers $2y'_{[(M+1)/2]}+1>\do>2y'_2+1>2y'_1+1$ given by

$\{2p_t+\ps(t)+\s-t;t\in[1,\s],\k_t\ne\k_{\ps(t)}\}$ in case (i),

$\{2p_t+\ps(t)+\s-t;t\in[1,\s],\k_t=\k_{\ps(t)}\}$ in case (ii),

$\{2p_t+\ps(t)+\s-t+1;t\in[1,\s],\k_t=\k_{\ps(t)}\}\sqc\{1\}$ in case (iii).
\nl
Thus, the sets $(\{y'_{[(M+1)/2]}>\do>y'_2>y'_1\},\{y_{[M/2]}>\do>y_2>y_1\})$ are given by
$$\align&(\{p_t+\t+(\ps(t)-t-1)/2;t\in[1,\s],\k_t\ne\k_{\ps(t)}\},\\&
\{p_t+\t+(\ps(t)-t)/2;t\in[1,\s],\k_t=\k_{\ps(t)}\})\\&=(\{p_t+\t+(-t-2)/2;t\in[1,\s],\ps(t)=-1,\k_t=0 \}\\&
\sqc\{p_t+\t+(-t-1)/2;t\in[1,\s],\ps(t)=0,\k_t=1 \},\\&
\{p_t+\t+(1-t)/2;t\in[1,\s],\ps(t)=1,\k_t=1\}\\&\sqc\{p_t+\t+(-t)/2;t\in[1,\s],\ps(t)=0,\k_t=0\})\\&
=(\{p_t+\t+(-t-2)/2;t\in[1,\s],\ps(t)=-1,\k_t=0\}\\&\sqc\{p_{t'}+\t+(-t'-2)/2;t'\in[1,\s],\ps(t')=0,\k_{t'}=0\},
\\&\{p_t+\t+(1-t)/2;t\in[1,\s],\ps(t)=1,\k_t=1\}\\&\sqc\{p_{t'}+\t+(1-t')/2;t'\in[1,\s],\ps(t')=0,\k_{t'}=1\})\\&
=(\{p_t+\t+(-t-2)/2;t\in[1,\s],\k_t=0\}\\&\sqc\{p_t+\t+(1-t)/2;t\in[1,\s],\k_t=1\})\endalign$$
in case (i),
$$\align&(\{p_t+\t+(\ps(t)-t)/2;t\in[1,\s],\k_t=\k_{\ps(t)}\},\\&\{p_t+\t+(\ps(t)+1-t)/2;t\in[1,\s],
\k_t\ne\k_{\ps(t)}\})
\\&=(\{p_t+\t+(1-t)/2;t\in[1,\s],\ps(t)=1,\k_t=1\}\\&\sqc\{p_t+\t-t/2;t\in[1,\s],\ps(t)=0,\k_t=0\},\\&
\{p_t+\t-t/2;t\in[1,\s],\ps(t)=-1,\k_t=0\}\\&\sqc\{p_t+\t+(1-t)/2;t\in[1,\s],\ps(t)=0,\k_t=1\})\\&
=(\{p_t+\t+(1-t)/2;t\in[1,\s],\ps(t)=1,\k_t=1\}\\&\sqc\{p_{t'}+\t+(1-t')/2;t'\in[1,\s],\ps(t')=0,\k_{t'}=1\},\\&
\{p_t+\t-t/2;t\in[1,\s],\ps(t)=-1,\k_t=0\}\\&\sqc\{p_{t'}+\t-t'/2;t'\in[1,\s],\ps(t')=0,\k_{t'}=0\})\\&
=(\{p_t+\t+(1-t)/2;t\in[1,\s],\k_t=1\},\\&\{p_t+\t-t/2;t\in[1,\s],\k_t=0\})\endalign$$
in case (ii),
$$\align&(\{p_t+\t+(\ps(t)-t)/2;t\in[1,\s],\k_t=\k_{\ps(t)}\}\sqc\{0\},\\&\{p_t+\t+(\ps(t)+1-t)/2;t\in[1,\s],
\k_t\ne\k_{\ps(t)}\})=\\&(\{p_t+\t+(1-t)/2;t\in[1,\s],\ps(t)=1,\k_t=1\}\sqc\{0\}\\&\sqc
\{p_t+\t-t/2;t\in[1,\s],\ps(t)=0,\k_t=0\}\sqc\{0\},\\&
\{p_t+\t-t/2;t\in[1,\s],\ps(t)=-1,\k_t=0\}\\&\sqc\{p_t+\t+(1-t)/2;t\in[1,\s],\ps(t)=0,\k_t=1\})\\&
=(\{p_t+\t+(1-t)/2;t\in[1,\s],\ps(t)=1,\k_t=1\}\sqc\{0\}\\&\sqc
\{p_{t'}+\t+(1-t')/2;t'\in[1,\s],\ps(t')=0,\k_{t'}=1\}\sqc\{0\},\\&
\{p_t+\t-t/2;t\in[1,\s],\ps(t)=-1,\k_t=0\}\\&\sqc\{p_{t'}+\t-t'/2;t\in[1,\s],\ps(t')=0,\k_{t'}=0\})\\&
=(\{p_t+\t+(1-t)/2;t\in[1,\s],\k_t=1\}\sqc\{0\},\\&\{p_t+\t-t/2;t\in[1,\s],\k_t=0\}\})\endalign$$
in case (iii). (We have used that, if $\k_t=0$, $t<\s$ and $\ps(t)=0$, then $ps(t+1)=\ps(t),p_{t+1}=p_t$; if
$\k_t=1$ and $\ps(t)=0$ then $\ps(t-1)=\ps(t),p_{t-1}=p_t$.)

In case (i) the multisets 
$$(\{y'_\t-\t+1\ge\do\ge y'_2-1\ge y'_1\},\{y_\t-\t+1\ge\do>y_2-1\ge y_1\})$$ 
are given by
$$(\{p_2-1\ge p_4-1\ge\do\ge p_\s-1\},\{p_1+1\ge p_3+1\ge\do\ge p_{\s-1}+1\}).$$
In case (ii) the multisets 
$$(\{y'_{\t+1}-\t\ge\do\ge y'_2-1\ge y'_1\},\{y_\t-\t+1\ge\do\ge y_2-1\ge y_1\})$$ 
are given by 
$$(\{p_1\ge p_3\ge\do\ge p_\s\},\{p_2\ge p_4\ge\do\ge p_{\s-1}\}).$$
In case (iii) the multisets 
$$(\{y'_{\t+1}-\t\ge\do\ge y'_2-1\ge y'_1\},\{y_\t-\t+1\ge\do\ge y_2-1\ge y_1\})$$ 
are given by 
$$(\{p_1\ge p_3\ge\do\ge p_{\s-1}\ge0\},\{p_2\ge p_4\ge\do\ge p_\s\}).$$
Now (a) follows from \cite{\IC, \S13}. Using (a) and 3.2(i),(h), we see that
$$\align&\r_\g=j_{S_{\bp_k}\T\do\T S_{\bp_2}\T W'_{\bp_1/2}\T W'_{\bp_1/2}}^{W'_n}(\sg)\text{ in case (i)};\\&
\r_\g=j_{S_{\bp_k}\T\do\T S_{\bp_2}\T W_{(\bp_1-\k_\s)/2}\T W'_{(\bp_1+\k_\s)/2}}^{W_n}(\sg)\\&
\text{ in case (ii),(iii)}.\tag b\endalign$$

\subhead 3.5\endsubhead
We prove Theorem 0.6(ii). If $G$ is of type $A$ the result is immediate. If $G$ is of classical type other than
$A$, the result follows from 3.3(b) and 3.4(b). If $G$ is of exceptional type, the result follows from the tables
in 2.3 using the data on $j$-induction available in the CHEVIE package \cite{\GH}. This completes the proof of 
Theorem 0.6.

\subhead 3.6\endsubhead
Let $w_{p_*}$ be the conjugacy class in $W_n$ (as in 3.1) associated to $p_*=(p_1\ge p_2\ge\do\ge p_\s)$ (where 
$p_1+p_2+\do+p_\s=n)$ as in \cite{\WEUN, 1,6}. We can view $w_{p_*}$ as an element of $\WW$ in the cases where 
$G$ is as in 1.2 with either $Q=0,\k=0$ or with $Q\ne0,\k=1$. In both cases $w_{p_*}$ has minimal length in its 
conjugacy class $C$ which is in $\uWW_{el}$. Let $\g_C$, $\g'_C$ be the corresponding $C$-small unipotent classes
(one is in a symplectic group, one is in an odd orthogonal group). From 3.3(b), 3.4(b), we see, using 3.2(g),(h),
that

(a) {\it the Springer representations $\r_{\g_C}$, $\r_{\g'_C}$ are the same.}
\nl
We see that the map $C\m\r_{\g_C}$ from $\uWW_{el}$ to $\Irr(\WW)$ depends only on the Weyl group $\WW$, not on
the underlying root system.


\widestnumber\key{GP}
\Refs
\ref\key\DL\by P.Deligne and G.Lusztig\paper Representations of reductive groups over finite fields\jour 
Ann.Math.\vol103\yr1976\pages103-161\endref
\ref\key\GH\by M.Geck, G.Hiss, F.L\"ubeck, G.Malle and G.Pfeiffer\paper A system for computing and processing
generic character tables for finite groups of Lie type, Weyl groups and Hecke algebras\jour Appl. Algebra Engrg.
Comm. Comput.\vol7\yr1996\pages1175-210\endref
\ref\key\GP\by M.Geck and G.Pfeiffer\book Characters of finite Coxeter groups and Iwahori-Hecke \lb algebras\publ 
Clarendon Press Oxford\yr2000\endref
\ref\key\OR\by G.Lusztig\book Characters of reductive groups over a finite field\bookinfo Ann.Math. Studies 107
\publ Princeton U.Press\yr1984\endref
\ref\key\IC\by G.Lusztig\paper Intersection cohomology complexes on a reductive group\jour Inv.Math.\vol75\yr1984
\pages205-272\endref
\ref\key\CSI\by G.Lusztig\paper Character sheaves, I\jour Adv.Math.\vol56\yr1985\pages193-237\endref
\ref\key\CSIII\by G.Lusztig\paper Character sheaves, III\jour Adv.Math.\vol57\yr1985\pages266-315\endref
\ref\key\SPE\by G.Lusztig\paper Unipotent classes and special Weyl group representations\jour J.Algebra\vol321
\yr2009\pages3418-3449\endref
\ref\key\WEUN\by G.Lusztig\paper From conjugacy classes in the Weyl group to unipotent classes\jour 
arxiv:1003.0412\endref
\ref\key\SPS\by N.Spaltenstein\paper On the generalized Springer correspondence for exceptional groups\inbook
Algebraic groups and related topics, Adv.Stud.Pure Math.6\publ North Holland and Kinokuniya\yr1985\pages317-338
\endref
\endRefs
\enddocument